\documentclass{amsart}
\usepackage{amssymb}
\usepackage{color}
\usepackage{hyperref}
\hypersetup{
    linktocpage=true,
    colorlinks=true,
    linkcolor=blue,
    filecolor=magenta,      
    urlcolor=cyan,
    pdftitle={Isometries between quantum symmetric spaces},
    bookmarks=true,
   citecolor=red    }

\newtheorem{thm}{Theorem}[section]
\newtheorem{lem}[thm]{Lemma}
\newtheorem{prop}[thm]{Proposition}
\newtheorem{cor}[thm]{Corollary}

\theoremstyle{definition}

\theoremstyle{remark}
\newtheorem{rem}[thm]{Remark}

\numberwithin{equation}{section}

\newcommand{\tma}{S(\mathcal{A},\tau)} 

\newcommand{\A}{\mathcal{A}}
\newcommand{\B}{\mathcal{B}}

\newcommand{\pa}{\mathcal{P}(\mathcal{A})}    %projections on \A%
\newcommand{\pb}{\mathcal{P}(\mathcal{B})}    %projections on \B%
\newcommand{\paf}{\mathcal{P}(\mathcal{A})_f}    %trace-finite projections on \A%

\newcommand{\dfs}[2]{d\left({#1}\right)({#2})}   

\newcommand{\svf}[1]{\mu_{{#1}}}  

\newcommand{\svft}[2]{\mu_{{#1}}({#2})}  

\newcommand{\subm}{\prec \prec}

\newcommand{\norm}[1]{\bigl\Vert{#1} \bigr\Vert}

\newcommand{\limm}[1]{\underset{#1}{\lim}}

\newcommand{\ra}{\rightarrow}

\newcommand{\id}{\mathbf{1}} 

\newcommand{\net}[1]{\{{#1}_{\lambda}\}_{\lambda \in \Lambda}}

\newcommand{\sotc}{\overset{SOT}{\rightarrow}}  

\newcommand{\seq}[1]{({#1}_n)_{n=1} ^{\infty}}
   
\newcommand{\summ}[2]{\underset{#1}{\overset{#2}{\sum}}}

\newcommand{\rax}[1]{\overset{#1}{\rightarrow}}

\newcommand{\G}{\mathcal{G}(\mathcal{A})}

\newcommand{\mtop}{\mathcal{T}_m}

\newcommand{\xmx}[2]{S({#2},{#1})} %()-measurable operators on ()%

\newcommand{\D}{\mathcal{D}}
\newcommand{\jn}{\vee}
\newcommand{\n}{\mathbb{N}}
\newcommand{\br}{\mathbb{R}}
\newcommand{\nmb}{S(\mathcal{B},\nu)} %tau-measurable operators%
\newcommand{\sotlimx}[1]{\text{SOT}\underset{#1}{\lim}}
\newcommand{\lmtop}{\mathcal{T}_{lm}}

\newcommand{\lmtopc}{\overset{\mathcal{T}_{lm}}{{\rightarrow}}}
\newcommand{\mtopc}{\overset{\mathcal{T}_{m}}{{\rightarrow}}}
\newcommand{\px}[1]{\mathcal{P}(\mathcal{#1})}    %projections on ?%
\newcommand{\uparrowx}[1]{\uparrow_{#1}}
\newcommand{\bhx}[1]{\mathcal{B}(#1)}

\newcommand{\unionx}[2]{\underset{#1}{\overset{#2}{\cup}}}

\begin{document}

\title[Isometries between non-commutative symmetric spaces]{Isometries between non-commutative symmetric spaces associated with semi-finite von Neumann algebras}


%    author one information
\author{Pierre de Jager}
\address{DST-NRF CoE in Math. and Stat. Sci\\ Unit for BMI\\ Internal Box 209, School of Comp., Stat., $\&$ Math. Sci.\\
NWU, PVT. BAG X6001, 2520 Potchefstroom\\ South Africa}
%\curraddr{}
\email{28190459@nwu.ac.za}
\thanks{The first author would like to thank the NRF for funding towards this project in the form of scarce skills and grantholder-linked bursaries}

%    author two information
\author{Jurie Conradie}
\address{Department of Mathematics, University of Cape Town, Cape Town, South Africa}
%\curraddr{}
\email{jurie.conradie@uct.ac.za}
%\thanks{}

%    \subjclass is required by all journals except JAG.

\date{\today}

%\dedicatory{}

%    The "communicated by" line appears only in PROC and JAG.
%\commby{}

%    Abstract is required.
\begin{abstract}
In this article we show that positive surjective isometries between symmetric spaces associated with semi-finite von Neumann algebras are projection disjointness preserving if they are finiteness preserving. This is subsequently used to obtain a structural description of such isometries. Furthermore, it is shown that if the initial symmetric space is a strongly symmetric space with absolutely continuous norm, then a similar structural description can be obtained  without requiring positivity of the isometry.
\end{abstract}

\maketitle

\section{Introduction}

The form of isometries between $L^p$-spaces was first described by Banach (in the case of finite measure spaces (\cite{key-Ban32})) and Lamperti (for $\sigma$-finite  measure spaces (\cite{key-Lam58})). In the proofs of these results  essential use is made of the fact that isometries map functions with disjoint support to functions with disjoint support. Representations of isometries between more general symmetric function spaces were obtained by Zaidenberg (\cite{key-Z}). We will define symmetric spaces below, but mention that well-known examples of such spaces include the $L^p$, Orlicz and Lorentz function spaces. A detailed account of results on isommetries in the commutative settings and the techniques used in the proofs can be found in \cite{key-F1}.

Non-commutative symmetric spaces  are Banach spaces of closed, densely-defined operators affiliated with a von Neumann algebra.  In the special case where the underlying von Neumann algebra is commutative, and hence isometrically isomorphic to an $L^\infty$ space over some localizable measure space, we obtain the commutative (classical) symmetric function spaces.  In the more general non-commutative (quantum) setting, isometries of  $L^p$-spaces associated with a semi-finite von Neumann algebra equipped with a faithful, normal semi-finite trace have been characterized by Yeadon (\cite{key-Y81}), but the description of isometries between  more general symmetric spaces have typically been limited to the finite trace setting or particular examples of semi-finite von Neumann algebras. In particular, structural descriptions for surjective isometries between Lorentz spaces (\cite{key-Chilin89}), positive surjective isometries between a symmetric space and a fully symmetric space (\cite{key-Chilin89}), and positive (not necessarily surjective) isometries between a symmetric space and a fully symmetric space with $K$-strictly monotone norm (\cite{key-Suk18}) have been obtained in the setting where the von Neumann algebra is equipped with a finite trace. Furthermore, surjective isometries on a separable symmetric space have been characterized (\cite{key-Suk96}) under the assumption that the underlying von Neumann algebra is an AFD (almost finite-dimensional) factor of type $II_1$ or $II_\infty$. In this paper we complement these results by considering surjective isometries between (general) symmetric spaces associated with (general) semi-finite von Neumann algebras.

% and these spaces generalize commutative and non-commutative $L^p$, Lorentz, Orlicz and Banach function spaces.

The technique we will employ is to analyze and utilize disjointness preserving properties of isometries. The motivation is as follows. Every von Neumann algebra is generated by its lattice of projections and therefore it is unsurprising that any isometric isomorphism between von Neumann algebras has to be implemented by a map that preserves this lattice structure, namely a Jordan $*$-isomorphism, possibly multiplied by a unitary operator (\cite{key-K51}). Furthermore, one would anticipate that there would be a relationship between the isometries of symmetric spaces associated with semi-finite von Neumann algebras and the isometries of the underlying von Neumann algebras. In describing the structure of an isometry between  symmetric spaces it is therefore natural to use the isometry to initially define a map on projections. In order to ensure that this map preserves the projection lattice structure and can be extended in a well-defined and linear manner, this map should preserve orthogonality of projections. In the setting of commutative and non-commutative $L^p$-spaces, for example, this can be achieved by showing that the isometry is disjointness preserving (\cite{key-Lam58} and \cite{key-Y81}, respectively). More recently it has been shown (\cite{key-Suk18}) that a positive isometry $T:E\ra F$ between symmetric spaces associated with semi-finite von Neumann algebras is disjointness preserving provided $F$ is contained in $L^0(\tau)$ and $F$ has $K$-strictly monotone norm (definitions to follow). This result is then used to describe the structure of a positive isometry $T:E \ra F$, where $E$ is a symmetric space on a trace-finite von Neumann algebra and $F$ is a fully symmetric space with $K$-strictly monotone norm on a trace-finite von Neumann algebra. In this paper we define a weaker notion of projection disjointness preserving maps, identify positive isometries satisfying this condition and show that even in the semi-finite setting, this weaker notion is sufficient to describe the structure of such isometries. 

The structure of the paper is as follows. In $\S \ref{S3}$ we obtain a local representation of positive surjective isometries, which enables us to show that these isometries are projection disjointness preserving. We then investigate projection disjointness preserving isometries in $\S \ref{S4}$ and show that even if these are not necessarily positive nor surjective we can describe their structure on an ideal contained in the intersection of the von Neumann algebra and the symmetric space. In order to obtain a global representation we consider isometries with more structure for the remainder of $\S \ref{S4}$. In $\S \ref{S5}$ we show that we can also obtain a global representation of projection disjointness preserving isometries with fewer assumptions on their structure if the initial symmetric space has slightly more structure.

Most of results in this paper will be proved under the assumption that the isometry under consideration is what we will call \emph{finiteness preserving}.  It will be shown in a subsequent paper (\cite{key-dJ19b}) that surjective isometries between Lorentz spaces associated with semi-finite von Neumann algebras satisfy this condition (and are also projection disjointness preserving). Furthermore, this condition is trivially satisfied if the final von Neumann algebra is equipped with a finite trace.

\section{Preliminaries}

Throughout this paper, unless indicated otherwise, we will use $\A\subseteq B(H)$ and $\B\subseteq B(K)$ to denote semi-finite von Neumann algebras, where $B(H)$ and $B(K)$ are the spaces of all bounded linear operators on Hilbert spaces $H$ and $K$, respectively. Let $\tau$ and $\nu$ denote distinguished faithful normal semi-finite traces on $\A$ and $\B$, respectively. The lattice of all projections in $\A$ will be denoted $\pa$ and the sublattice of projections with finite trace will be denoted $\paf$. We will use $\id$ to denote the identity of $\A$. The set of all finite linear combinations of mutually orthogonal projections in $\pa$ (alternatively $\paf$) will be denoted $\mathcal{G}(\mathcal{A})$ (respectively $\mathcal{G}(\mathcal{A})_f$).  Convergence in $\A$ with respect to the  operator norm topology, the strong operator topology (SOT) and the weak operator topology (WOT) will be denoted by respectively $\rax{\A},\; \rax{SOT}$ and $\rax{WOT}$. A linear map $\Phi:\A\to \B$ is called a \textit{Jordan homomorphism} if $\Phi(yx + xy) = \Phi(y)\Phi(x) + \Phi(x)\Phi(y)$ for all $x, y \in \A$. If, in addition, $\Phi(x^\ast) = \Phi(x)^\ast$ for all $x\in \A$, then $\Phi$ is called a \textit{Jordan $\ast$-homomorphism}. Further details regarding von Neumann algebras and Jordan homomorphisms may be found in \cite{key-K1}.

A closed  operator $x$ with domain $\D(x)$ dense in $H$ is \textit{affiliated} with $\A$ if $u^*xu=x$ for all unitary operators $u$ in the commutant $\A^\prime$ of $\A$. A closed densely defined self-adjoint operator $x$ with  spectral measure $e^x$ is affiliated to $\A$ iff $e^x(B)\in \pa$ for every Borel subset $B$ of $\mathbb{R}$. For such an operator we will write $x=\int_{-\infty}^{\infty} \lambda de^x_\lambda$ if $\{e^{x}_\lambda\}_\lambda$ is the unique resolution of the identity such that $x\eta = \int_{-n}^{n} \lambda de^x_\lambda\eta$ for each $\eta \in f_n(H)$ and all $n$, and $\unionx{n=1}{\infty}f_n(H)$  is a core for $x$, where $f_n:=e^x_{n}-e^{x}_{-n}$ (see \cite[Theorem 5.6.12]{key-K1}). If $x: \D(x)\to H$ is a closed and densely defined operator, then the projection onto the kernel of $x$ will be denoted by $n(x)$, the projection onto closure of the range of $x$ by $r(x)$, and the \textit{support projection} $\id - n(x)$ by $s(x)$. It follows that $x=r(x) x= x s(x)$, and if $x=x^\ast$, then $r(x)=s(x)$ and $x=s(x) x =x s(x)$. If $x$ is affiliated with $\A$, all three these projections are in $\A$. A closed, densely defined operator $x$ affiliated to $\A$ is called \textit{$\tau$-measurable} if there is a sequence $(p_n)$ in $\pa$ such that $p_n \uparrow \id$, $p_n(H)\subseteq \D(x)$ and $\id - p_n\in\paf$ for every $n$. It is known that if $x=u |x|$  is the polar decomposition of $x$, then $x$ is $\tau$-measurable if and only if it is affiliated to $\A$ and there is a $\lambda>0$ such that $\tau(e^{|x|}(\lambda,\infty))<\infty$. A vector subspace $\D \subseteq H$ is is called $\tau$-dense  if there exists a sequence $(p_n)$ in $\pa$ such that $p_n(H) \subseteq \D$ for all $n$, $p_n \uparrow 1$ and $\tau(\id - p_n) < \infty$ for all $n$. Clearly a closed densely defined  operator $x$ affiliated to $\A$ is $\tau$-measurable if and only its domain $\D(x)$ is $\tau$-dense. The set of all $\tau$-measurable operators affiliated with $\A$ will be denoted $\tma$ or $S(\A)$. It becomes a $\ast$-algebra when sums and products are defined as the closures of respectively the algebraic sum and algebraic product. For $x\in \tma$ we write $x\geq 0$ if  $\langle x\xi,\xi\rangle\geq 0$ for all $\xi$ in the domain of $x$ (where $\langle\cdot,\cdot\rangle$ denotes the inner product on $H$), and we put $\tma^+=\{x\in\tma: x\geq 0\}$. The cone $\tma^+$ defines a partial order on the self-adjoint elements of $\tma$. If $\mathcal{H}$ is any collection of $\tau$-measurable operators, then we will write $\mathcal{H}^{sa}=\{x\in\mathcal{H}:x=x^\ast\}$ and $\mathcal{H}^+=\{x\in\mathcal{H}:x\geq 0\}$. Note that $\A$ is an \emph{absolutely solid subspace} of $\tma$, i.e. if $x\in \tma$ and $y\in \A$ with $|x|\leq |y|$, then $x\in \A$.

For $\epsilon,\delta>0$, define $N(\epsilon,\delta):=\{x\in \tma: \tau(e^{|x|}(\epsilon,\infty))\leq \delta\}$. The collection $\{N(\epsilon,\delta):\epsilon,\delta>0\}$ defines a neighbourhood base for a vector space topology $\mtop$ on $\tma$. This topology is called the \emph{measure topology} and with respect to this topology $\tma$ is a complete metrisable topological $*$-algebra. We will repeatedly use the fact that multiplication is jointly continuous in the measure topology. Another important vector space topology on $\tma$ is the \emph{local measure topology}, denoted $\lmtop$, which has a neighbourhood base consisting of the collection of sets of the form $N(\epsilon, \delta, p):=\{x\in\tma: pxp\in N(\epsilon,\delta)\}$, where $\epsilon,\delta>0$ and $p\in \paf$. Multiplication is separately, but not jointly continuous with respect to the local measure topology, that is $x_\alpha y \lmtopc xy$ and $yx_\alpha \lmtopc yx$ whenever $y\in \tma$ and $\{x_\alpha\}_\alpha$ is a net in $\tma$ with $x_\alpha \lmtopc x\in \tma$. 

If $(x_\lambda)_{\lambda \in \Lambda}$ is an increasing net in $\tma$ and $x=\sup\{x_\lambda:\lambda \in \Lambda\}\in\tma$, we write $x_\lambda \uparrow x$. In the case of a decreasing net $(x_\lambda)_{\lambda \in \Lambda}$ with infimum $0$ we write $x_\lambda\downarrow 0$. If $\mathcal{H}\subseteq \tma$ and $T:\mathcal{H} \ra \nmb$ is a linear map such that $T(x_\lambda)\uparrow T(x)$ whenever $\{x_\lambda\}_{\lambda \in \Lambda}$ is a net in $\mathcal{H}^{sa}$  such that $x_\lambda \uparrow x\in \mathcal{H}^{sa}$, then $T$ will be called \emph{normal} (on $\mathcal{H}$). If $E$ is a linear subspace of $\tma$, a linear map $T:E\ra \nmb$ will be called \emph{finiteness preserving} if $\nu(s(T(p))<\infty$ whenever $p\in \paf$. For background and further details regarding trace-measurable operators the interested reader is referred to \cite{key-Dodds14} and \cite{key-Terp1}. 
 
For $x\in \tma$, the \emph{distribution function} of $|x|$ is defined as $\dfs{|x|}{s}:=\tau\left(e^{|x|}(s,\infty)\right)$, for $s\geq 0$. The \emph{singular value function} of $x$, denoted $\svf{x}$, is defined to be the right continuous inverse of the  distribution function of $|x|$, namely \[\svft{x}{t}=\inf\{s\geq 0: \dfs{|x|}{s}\leq t\} \qquad t\geq 0. \]  If $x,y \in S(\mathcal{A},\tau)$, then we will say that $x$ is \emph{submajorized} by $y$ and write $x \prec\prec y$ if $\int_0^t \svft{x}{s}ds\leq \int_0^t \svft{y}{s}ds$ for all $t> 0$.  Let $S^0(\A,\tau)$ denote the ideal of $\tau$-compact operators, which is defined as the set of all $\tau$-measurable operators $x$ for which $\limm{t\ra \infty} \svft{x}{t}=0$. 

A linear subspace $E\subseteq S(\mathcal{A},\tau)$, equipped with a norm $\norm{\cdot}_E$, is called a \emph{symmetric space} if  $E$ is a Banach space  and $x\in E$ with $\norm{x}_E\leq \norm{y}_E$, whenever $y\in E$ and $x\in S(\mathcal{A},\tau)$ with $\svf{x}\leq \svf{y}$. In this case we also have that $uxv\in E$ and $\norm{uxv}_E\leq \norm{u}_{\mathcal{A}}\norm{v}_{\mathcal{A}}\norm{x}_E$ for all $x\in E,u,v\in \mathcal{A}$. Furthermore, $\norm{x}_E=\norm{x^\ast}_E=\norm{|x|}_E$ for all $x\in E$, and  $\norm{x}_E\leq \norm{y}_E$ whenever  $x,y\in E$ with $|x|\leq |y|$. A symmetric space is an absolutely solid subspace of $\tma$.  A symmetric space $E\subseteq S(\mathcal{A},\tau)$ is called \emph{strongly symmetric} if its norm has the additional property that $\norm{x}_E\leq \norm{y}_E$, whenever $x,y\in E$ satisfy $x \prec\prec y$.  If $E$ is a symmetric space and it follows from $x\in S(\mathcal{A},\tau)$, $y \in E$ and $x \prec\prec y$ that $x\in E$ and $\norm{x}_E \leq \norm{y}_E$, then $E$ is called a \emph{fully symmetric space}. Let $E \subseteq S(\mathcal{A},\tau)$ be a symmetric space. Convergence in $E$ with respect to the norm of $E$  will be denoted by $\rax{E}$. The \emph{carrier projection} $c_E$ of $E$ is defined to be the supremum of all projections in $\mathcal{A}$ that are also in $E$. If $c_E=\id$, then $E$ is continuously embedded in $S(\mathcal{A},\tau)$ equipped with the measure topology $\mathcal{T}_m$. We will assume throughout this paper that $c_E=\id$. The norm $\norm{\cdot}_E$ on a symmetric space $E$ is called order continuous if $\norm{x_\lambda}\downarrow 0$ whenever $x_\lambda \downarrow 0$ in $E$. If this is the case, $\mathcal{F}(\tau):=\{x\in\A: s(x)\in \paf \}$ is norm dense in $E$, and it can be shown, using the spectral theorem, that for every $x\in \A^{sa}$, there is a sequence $(x_n)_{n=1}^{\infty}$ in $\G_f$ such that $x_n\rax{E} x$. If $E$ is a strongly symmetric space, then it can be shown (\cite[Proposition 6.12]{key-DdP12}) that $E$ has order continuous norm if and only if it has absolutely continuous norm, that is $\norm{p_nxp_n}_E \ra 0$ for every sequence $\seq{p}$ in $\pa$ satisfying $p_n \downarrow 0$ and every $x\in E$. 

If $\A=L^\infty(0,\infty)$ is the abelian semi-finite von Neumann algebra of all essentially bounded Lebesgue measurable functions on $(0,\infty)$ and the trace $\tau$ is given by integration with respect to Lebesgue measure, then $\tma=\mathcal{S}(0,\infty)$ is the space of all Lebesgue measurable functions on $(0,\infty)$ that are bounded except possibly on a set of finite measure. In this case the singular value function $\svf{x}$ corresponds to the decreasing rearrangement $f^\ast$ of a measurable function $f$. It follows from \cite[Corollaries 2.6 and 2.7]{key-Dodds92} that if $(\mathcal{A},\tau)$ is a semi-finite von Neumann algebra and $E(0,\infty)\subseteq \mathcal{S}(0,\infty)$ is a fully symmetric space, then the set $E(\A):=\{x\in S(\mathcal{A},\tau):\svf{x}\in E(0,\infty)\}$ is a fully symmetric space, when equipped with the norm $\norm{x}_{E(\A)}=\norm{\svf{x}}_{E(0,\infty)}$ for $x\in E(\A)$. Furthermore, similar results hold for symmetric spaces and strongly symmetric spaces (see \cite{key-Kalton08} and \cite{key-Dodds14}).

The following easily verifiable result will be used repeatedly and details conditions under which convergence in a von Neumann algebra yields convergence in an associated symmetric space.

\begin{prop} \label{P8 05/08/14} 
Suppose $E\subseteq \tma$ is a symmetric space. If $\seq{x}$ is a sequence in $E \cap \A$ is such that $x_n \rax{\A} x\in E\cap \A$ and either $s(x),s(x_n)\leq p$ or $r(x),r(x_n)\leq p$ for all $n \in \n^+$ and for some $p\in \paf$, then $x_n \rax{E} x$.
\end{prop}

Since any symmetric space $E\subseteq \tma$ is continuously embedded in $\tma$ equipped with the measure topology (\cite[Proposition 20]{key-Dodds14}), we obtain the following corollary.  

\begin{cor} \label{CP8 05/08/14} Suppose $(\A,\tau)$ and $(\B,\nu)$ are semi-finite von Neumann algebras and $E\subseteq S(\mathcal{A},\tau)$ and $F\subseteq S(\mathcal{B},\nu)$  are symmetrically normed spaces. If $U:E\ra F$ is a continuous map with respect to the norms on $E$ and $F$, then $U(x_n)\rax{\mtop} U(x)$, whenever $\seq{x}$ is a sequence in $\mathcal{F}(\tau)$ such that $x_n\rax{\A} x\in \mathcal{F}(\tau)$ and $s(x_n)\leq s(x)$ or  $r(x_n)\leq r(x)$ for all $n\in \mathbb{N}^+$).
\end{cor}

In \cite{key-Suk18}, a linear map $U:E\ra F$ between symmetric spaces is called \emph{disjointness preserving} if $U(x)U(y)=0$ whenever $x,y\in E^+$ with $xy=0$. For the purposes of this paper we introduce a slightly weaker notion. We will call a linear map $U:E\subseteq \tma \ra \nmb$ \emph{projection disjointness preserving} if $U(p)^*U(q)=U(p)U(q)^*=0$, whenever $p,q\in \paf$ with $pq=0$. It is clear that a positive map will be projection disjointness preserving whenever it is disjointness preserving. We provide sufficient conditions for the converse to hold.

\begin{prop}\label{P2 03/12/18}
Suppose $E\subseteq \tma$ and $F\subseteq \nmb$ are symmetric spaces and $U:E \ra F$ is a bounded linear projection disjointness preserving map. If $E$ is strongly symmetric with absolutely continuous norm, or $F\subseteq S^0(\B, \nu)$ and $U$ is normal, then $U$ is disjointness preserving.
\begin{proof}
Suppose $x=\summ{i=1}{n}\alpha_i p_i,y=\summ{j=1}{m}\beta_j q_j\in \mathcal{G}(\mathcal{A})_f^+$ with $xy=0$. Then it is easily checked that $s(x)s(y)=0$, and for every $i,j$ we have that $p_iq_j=0$, since $p_i\leq s(x)$ and $q_j\leq s(y)$. Using the linearity and projection disjointness preserving nature of $U$ we therefore have that $U(x)U(y)=0$.

If $E$ has absolutely continuous norm and $x,y\in E^+$, then there exists $\seq{x},$ \\ $\seq{y}\subseteq \mathcal{G}(\A)_f^+$ such that $x_n \rax{E} x$ and $y_n \rax{E} y$. Therefore $U(x_n)\rax{F} U(x)$ and $U(y_n)\rax{F} U(y)$. By \cite[Proposition 20]{key-Dodds14}, this implies that $U(x_n)\mtopc U(x)$ and $U(y_n) \mtopc U(y)$ and so $U(x_n)U(y_n)\mtopc U(x)U(y)$, since multiplication is jointly continuous in the measure topology (\cite[p. 210]{key-Dodds14}). Furthermore, $s(x)x_n s(x) \rax{E} s(x)xs(x)=x$ and similarly $s(y)y_n s(y) \rax{E} y$. We can therefore assume without loss of generality that  $s(x_n)s(y_n)=0$ for every $n$ and thus $x_ny_n=0$ for every $n$. It follows that $U(x_n)U(y_n)=0$ for every $n$ and so $U(x)U(y)=0$.

If $F\subseteq S^0(\B,\nu)$ and $U$ is normal, then we first note that if $x,y\in \mathcal{F}(\tau)^+$, then there exists $\seq{x},\seq{y}\subseteq \mathcal{G}(\mathcal{A})_f$ such that $x_n \rax{\A} x$, $y_n\rax{\A} y$, $s(x_n)\leq s(x)$ and $s(y_n)\leq s(y)$ for every $n$. By Proposition \ref{P8 05/08/14}, $x_n \rax{E} x$ and $y_n \rax{E}y$. In the same way as before we can then show that $U(x)U(y)=0$. Finally, if $x,y \in E^+$, then there exists $\net{x}, \{y_\alpha\}_{\alpha \in \mathbb{A}} \subseteq \mathcal{F}(\tau)^+$ such that $x_\lambda \uparrow x$ and $y_\alpha \uparrow y$ (see \cite[p. 211]{key-Dodds14}). Since $U$ is normal we have that $U(x_\lambda) \uparrow U(x)$ and therefore $U(x_\lambda) \mtopc U(x)$, by \cite[Proposition 2(iv)]{key-Dodds14} (since $F\subseteq S^0(\B,\nu)$). Similarly, $U(y_\alpha) \mtopc U(y)$. Since $s(x_\lambda)\leq s(x)$ and $s(y_\alpha)\leq s(y)$ for each $\lambda$ and $\alpha$, we have that $U(x_\lambda)U(y_\alpha)=0$ for each $\lambda$ and $\alpha$ and so $U(x)U(y)=0$ as before. 
\end{proof}
\end{prop}

Further information about symmetric spaces may be found in \cite{key-Dodds14} and \cite{key-Pag}.

\section{The projection disjointness preserving property of positive surjective isometries} \label{S3}

In order to describe the structure of positive surjective isometries we will start by showing that under certain conditions such isometries are projection disjointness preserving. It is shown in \cite[Corollary 5]{key-Suk18} that if $T:E\ra F$ is a positive isometry, where $E\subseteq \tma$ is a symmetric space and $F\subseteq L^0(\B,\nu):=S^0(\B,\nu) \cap L^1+L^\infty(\B)$ is a symmetric space with $K$-strictly monotone norm, then $T$ is disjointness preserving. In this section we complement this result by showing that a finiteness preserving positive surjective isometry between arbitrary symmetric spaces is projection disjointness preserving. Suppose $E\subseteq \tma$ and $F\subseteq \nmb$ are symmetric spaces. We will start by showing that if $U:E\ra F$ is a positive surjective isometry, then $U$ is an order isomorphism and for each $p\in \paf$, $U$ maps $pEp$ into $s(U(p))Fs(U(p))$. Since we were not able to show that $U$ in fact maps $pEp$ onto $s(U(p))Fs(U(p))$ and we are not assuming full symmetry of $F$, we do not have access to \cite[Theorem 3.1]{key-Chilin89}, which would have enabled us to describe the structure of $U$ under the additional assumption that $U$ is finiteness preserving. Nevertheless, under this assumption we are able to adapt the technique employed in the proof of \cite[Theorem 3.1]{key-Chilin89} to prove a local representation of such isometries in the sense that for each $p\in \paf$ we will show that there exists a Jordan $*$-isomorphism $\Phi_p$ from $p\A p$ onto $s(U(p))\B s(U(p))$ such that $U(x)=U(p)\Phi_p(x)$ for all $x\in p\A p$. The projection disjointness preserving property of positive surjective isometries will then follow from this.

\begin{lem} \label{L1 C3.1} 
Suppose $E\subseteq \tma$ and $F\subseteq \nmb$ are symmetric spaces. If $U:E\ra F$ is a positive isometry, then  $z \geq 0$, whenever $z\in E$ and $U(z) \geq 0$. If in addition, $U$ is surjective, then $U$ is an order isomorphism and hence also normal.
\begin{proof}
The proof of the corresponding result in the setting where $F$ is a fully symmetric space and $\tau(\id),\nu(\id)<\infty$ (\cite[Lemma 3.2]{key-Chilin89}) requires only one significant adjustment to be generalized to spaces associated with arbitrary semi-finite von Neumann algebras. This proof uses the fact that if $\nu(\id)<\infty$, then $x-y\subm x+y$ whenever $x,y\in L^1(\B,\nu)$ (see \cite[Lemma 2.1]{key-Chilin89}). The full symmetry of $F$ is then used to show that $\norm{x-y}_F\leq \norm{x+y}_F$, if in addition $x,y\in F^+$.

To extend \cite[Lemma 3.2]{key-Chilin89} to the general semi-finite setting we note that it has recently been shown in \cite[Corollary 4]{key-Bik12} that, even in this more general setting, $\norm{x-y}_F\leq \norm{x+y}_F$ whenever $x,y\in F^+$ and $F$ is a normed solid space. Since symmetric spaces are normed solid spaces, we do not require the full symmetry assumption. Finally, it is easily checked that an order isomorphism is necessarily normal.
\end{proof}
\end{lem}

The following lemma will play an important role in obtaining a local representation of positive surjective isometries. 

\begin{lem}
Suppose $E\subseteq \tma$ and $F\subseteq \nmb$ are symmetric spaces and $U:E\ra F$ is a positive surjective isometry. If $p\in \paf$, then $U(pEp)\subseteq s(U(p))Fs(U(p))$.
\begin{proof} 
Since $U$ is positive we have that $s(U(p))=s(U(p)^*)=r(U(p))$. This implies that   $s(U(p))U(p)s(U(p))=U(p)$ and hence $U(p)\in s(U(p)) F s(U(p))$. If $q\in \pa$, then $0\leq pqp \leq p\id p$, by \cite[Proposition 1(iii)]{key-Dodds14} and so $0 \leq U(pqp) \leq U(p)$. This implies that $U(pqp)\in s(U(p))Fs(U(p))$. It follows that $U(p\mathcal{G}(\mathcal{A})p)\subseteq s(U(p))Fs(U(p))$. If $x\in pEp\cap \A \subseteq p\mathcal{F}(\tau)^+p$, then using the Spectral Theorem there exists $\seq{x}\subseteq \mathcal{G}(\mathcal{A})^+$ such that $x_n \rax{\A} x$ and $r(x_n)=s(x_n)\leq s(x)\leq p$ for each $n\in \mathbb{N}^+$. Then $x_n\in p\mathcal{G}(\mathcal{A})p$ for each $n$ and $U(x_n)\rax{F} U(x)$. Since $U(x_n)\in s(U(p))Fs(U(p))$ for each $n$ and it is easily checked that $s(U(p))Fs(U(p))$ is closed in $F$, we have that $U(x)\in  s(U(p))Fs(U(p))$.  Finally, if $x\in pE^+p$, then by \cite[Proposition 1(vii)]{key-Dodds14} there exists $\net{x}\subseteq \mathcal{F}(\tau)^+$ such that $x_\lambda \uparrow x$. Then $px_\lambda p \uparrow pxp =x$. It follows by Lemma \ref{L1 C3.1} that $U$ is normal and therefore $U(px_\lambda p)\uparrow U(x)$. It follows that $U(px_\lambda p) \lmtopc U(x)$, by \cite[Proposition 2(v)]{key-Dodds14}. Since $U(px_\lambda p)\in s(U(p))Fs(U(p))$ for each $\lambda$ and it is easily checked that $s(U(p))Fs(U(p))$ is closed in the local measure topology, we have that $U(x)\in s(U(p))Fs(U(p))$. 
\end{proof}
\end{lem}

Next we show how the techniques of \cite[$\S 3$]{key-Chilin89} may be adapted to obtain a local representation of positive surjective isometries. To facilitate this we mention a few aspects of reduced spaces (see \cite[p. 211, 212 and 215]{key-Dodds14}). For $p\in \paf$ and $x\in \tma$, let $x_{(p)}:=(pxp)\restriction_{p(H)}$, where $H$ denotes the Hilbert space on which $\A$ acts. It can be shown that $\{x_{(p)}:x\in \tma\}=\xmx{\tau_p}{\A_p}$, where $\A_p:=\{x_{(p)}:x\in \A\}$ and $\tau_p(x_{(p)}):=\tau(pxp)$ for every $x\in \A$.  Let $\phi_p$ denote the canonical map $x\mapsto x_{(p)}$ from $p\xmx{\tau}{\A}p$ onto $\xmx{\tau_p}{\A_p}$.  Note that $\phi_p$ is a $*$-isomorphism, $E_p$ is a symmetric space if $E$ is a symmetric space, and that the restrictions of $\phi_p$ to $p\A p$ and $pEp$ respectively are isometries onto the reduced spaces $\A_p=\{x_{(p)}:x\in \A\}$ and $E_p=\{x_{(p)}:x\in E\}$. Let $\psi_{p}$ denote the canonical map from $s(U(p))\xmx{\nu}{\B}s(U(p))$ onto $\xmx{\nu_{s(U(p))}}{\B_{s(U(p))}}$. We will make use of the fact that if $x\in p\tma^{sa} p$ and $f$ is a Borel measurable function on $\mathbb{R}$ that is bounded on compact sets, then $f (\phi_p(x))=\phi_p(f(x))$ and a similar relationship holds for elements in $s(U(p))\nmb^{sa}s(U(p))$ (this follows from an application of \cite[Proposition 2.9.2]{key-dP}).

\begin{prop} \label{PL4 C3.1}
Suppose $U:E\ra F$ is a positive surjective isometry. If $U$ is finiteness preserving, then for each $p\in \paf$,  there exists a Jordan $*$-isomorphism $\Phi_p$ from $p\A p$ onto $s(U(p))\B s(U(p))$ such that $U(x)=U(p)\Phi_p(x)$ for every $x\in p \A p$. Furthermore, $a_p:=U(p)$ commutes with every element in $s(U(p))S(\B) s(U(p))$. 
\begin{proof} 
For $p\in \paf$, let $a_p:=U(p)$. If we let $B_{s(U(p))}$ denote the reduced space corresponding to $s(U(p))\B s(U(p))$ and if we identify $a_p$ with the corresponding element in the reduced space $S(\B_{s(U(p))})\cong s(U(p))S(\B)s(U(p))$, then we have that $a_p$ is invertible in $S(\B_{s(U(p))})$ (this follows from the functional calculus for $a_p$ and noting that  $s(a_p)=s(U(p))$ is the identity of $\B_{s(U(p))}$ and has finite trace). We will use $a_p^{-1}$ denote the inverse of $a_p$ in $S(\B_{s(U(p))})$ (bearing in mind that $a_p$ need not be invertible in $S(\B)$). Working in these reduced spaces and using these identifications, we have that  $a_p^{-1}\geq 0$ and $a_p^{-1/2}=(a_p^{-1})^{1/2}=(a_p^{1/2})^{-1}$. In this setting, we let 
\[\Phi_p(x)=a_p^{-1/2}U(x)a_p^{-1/2} \qquad x\in \A_p. \]
Note that since $\A_p$ is trace-finite, $\A_p \subseteq E_p\cong pEp$ and so $U$ is defined on all of $\A_p$.  

It is easily checked that $\Phi_p$ is a positive unital map. To show that $\Phi_p$ maps $\A_p$ into $\B_{s(U(p))}$ note that if $y \in \A_p^+$, then $0\leq  y \leq \norm{y}_{\A_p}p$, by \cite[Proposition 4.2.3]{key-K1}. This implies that $0 \leq \Phi_p(y) \leq \norm{y}_{\A_p} \Phi_p(p)=\norm{y}_{\A_p} s(U(p))$  since $\Phi_p$ is positive, linear and unital. It follows that $\Phi_p(y)\in \B_{s(U(p))}$, since $\norm{y}_{\A_p}s(U(p)) \in \B_{s(U(p))}$ and $\B_{s(U(p))}$ is an absolutely solid subspace of $S(B_{s(U(p))})$. Since any element of $\A_p$ can be written as a linear combination of positive elements, we have that $\Phi_p(\A_p)\subseteq \B_{s(U(p))}$. Next we show that $\Phi_p$ is surjective. Let $b\in \B_{s(U(p))}^+$ and define $c=a_p^{1/2}ba_{p}^{1/2}$. Then
\begin{eqnarray}
&0 \leq c= a_p^{1/2}ba_{p}^{1/2} \leq a_p^{1/2}\norm{b}_{\B_{s(U(p))}}s(U(p))a_p^{1/2} = \norm{b}_{\B_{s(U(p))}} a_p.&   \label{e4.0 C3.1}
\end{eqnarray}
Since $F$ is symmetric, $F_{s(U(p))}$ is also symmetric. This, combined with (\ref{e4.0 C3.1}),  implies that $c\in F_{s(U(p))}$, since $\norm{b} a_p\in F_{s(U(p))}$. By Lemma \ref{L1 C3.1}, $U^{-1}$ is positive and therefore $0\leq U^{-1}(c)\leq \norm{b}_{\B_s(U(p))}p$, using (\ref{e4.0 C3.1}). It follows that $U^{-1}(c)\in p\A p$. Furthermore, it is easily checked that $\Phi_p(U^{-1}(c))=b$. It follows that $\Phi_p$ is surjective and for $y\in \B_{s(U(p))}$, $\Phi_p^{-1}(y)=U^{-1}(a_p^{1/2}ya_{p}^{1/2})$. Using this formula for the inverse of $\Phi_p$, \cite[Proposition 1(iii)]{key-Dodds14} and the positivity of $U^{-1}$, we see that $\Phi_p^{-1}$ is positive. We have shown that $\Phi_p$ is a unital order isomorphism of $\A_p$ onto $\B_{s(U(p))}$ and therefore $\Phi_p$ is a Jordan $*$-isomorphism, by \cite[Exercise 10.5.32]{key-K2}.
  
By definition of $\Phi_p$, we have that $\Phi_p(x)=a_p^{-1/2}U(x)a_p^{-1/2}$ and therefore $U(x) = a_p^{1/2}\Phi_p(x)a_p^{1/2}$. Essentially the same technique as the one employed in the proof of \cite[Lemma 3.5]{key-Chilin89} can be used to show that $a_p \in S(Z(\B_{s(U(p))}))$ (where $Z(\B_{s(U(p))})$ denotes the center of the von Neumann algebra $\B_{s(U(p))}$). It now follows that $a_p b=b a_p$ for every $b \in s(U(p))S(\B)s(U(p))$ and therefore $U(x)=a_p\Phi_p(x)$ for every $x\in p\B p$, by \cite[Proposition 2.2.22]{key-dP}.
\end{proof}
\end{prop}

\begin{cor} \label{CnnT3.1 Chilin89} 
Let $E\subseteq \tma$ and $F \subseteq \nmb$ be symmetric spaces and $U:E\ra F$ a positive surjective isometry. If $U$ is finiteness preserving (in particular if $\nu(\id)<\infty$), then $U$ is projection disjointness preserving. 
\begin{proof}
It follows from the previous result that if $p,q\in \paf$ with $pq=0$, then $p+q\in \paf$ and $U(p)U(q)=a_{p+q}^2\Phi_{p+q}(p)\Phi_{p+q}(q)=0$ (see \cite[Exercise 10.5.22(vii)]{key-K2}). 
\end{proof}
\end{cor}

\section{The structure of positive surjective isometries}\label{Positive surjective isometries} \label{S4}

Our aim in this section is to describe the structure of positive surjective isometries. We saw in the previous section that if, in addition, such an isometry is finiteness preserving, then it is projection disjointness preserving. We start by considering projection disjointness preserving isometries (that are not necessarily positive nor surjective). We show that the ideas of Yeadon's Theorem and the extension procedures developed in \cite{key-dJ18b} can be used to describe such isometries on $\mathcal{F}(\tau)$. More specifically we will show that if $V$ is a projection disjointness preserving isometry between symmetric spaces $E\subseteq \tma$ and $F\subseteq \nmb$, then letting $\Psi(p)=s(V(p))$ for $p\in \paf$ yields a projection mapping which can be extended to a positive linear map (still denoted $\Psi$) on $\mathcal{F}(\tau)$, which preserves squares of self-adjoint elements and therefore has many Jordan $*$-homomorphism-like properties (see \cite[Proposition 2.3]{key-dJ18b}). Furthermore, we will show that $V(x)=V(p)\Psi(x)=v_pb_p\Psi(x)$ for any $x\in \mathcal{F}(\tau)$ and $p\in \paf$ with $p\geq s(x)$, where $v_p$ and $b_p$ are respectively the partial isometry and positive operator occurring in the polar decomposition $V(p)=v_pb_p$. Attempts to extend $\Psi$ to all of $\A$ and use the $v_p$'s and $b_p$'s to construct single elements which can be used in a global representation of $V$ have proven to be problematic without further conditions on the symmetric spaces $E$ and $F$ or the isometry $V$. In this section we will show that the extension and representation can be achieved in the general setting of symmetric spaces if the isometry has more structure, and in the following section we will show how the extension and representation can be achieved if the isometry does not necessarily have all of this additional structure, provided the symmetric spaces have more structure. 

We will need the following extension result. 

\begin{thm}\label{nP extension2}\cite[Theorems 3.7 and 5.1]{key-dJ18b} \label{T2 04/07/18}\label{L2 03/07/18}
Suppose $\Phi:\paf \ra \pb$ is a map such that $\Phi(p+q)=\Phi(p)+\Phi(q)$ whenever $p,q\in \paf$ with $pq=0$. If there exists a linear map $U$ from $\mathcal{F}(\tau)$ into $\xmx{\nu}{\B}$ such that $\Phi(p)=s(U(p))$ for all $p\in \paf$, and which has the property that $U(x_n)\rax{\mtop} U(x)$ whenever $\seq{x}$ is a sequence in $\mathcal{F}(\tau)$ such that $x_n \rax{\A} x\in \mathcal{F}(\tau)$ and $s(x_n)\leq s(x)$ for all $n\in \mathbb{N}^+$, then $\Phi$ can be extended to a positive linear map (still denoted by $\Phi$) from $\mathcal{F}(\tau)$ into $\B$ such that $\norm{\Phi(x)}_\B\leq\norm{x}_\A$ and $\Phi(x^2)=\Phi(x)^2$ for all $x\in \mathcal{F}(\tau)^{sa}$. Suppose, in addition, that $U$ is positive and normal. 
\begin{enumerate}
\item If $x\in \mathcal{F}(\tau)^{sa}$ and $p\in \paf$ with $p\geq s(x)$, then  $\Phi(x)U(p)=U(x)=U(p)\Phi(x)=U(p)^{1/2}\Phi(x)U(p)^{1/2}$;
\item If $x\in \mathcal{F}(\tau)^{sa}$ and $p\in \paf$ with $p\geq s(x)$, then there exists a $w_p\in \nmb$ such that $U(p)^{1/2}w_p=\Phi(p)=w_pU(p)^{1/2}$ and $\Phi(x)=w_pU(x)w_p$;
\item  $\Phi$ can be extended to a normal Jordan $*$-homomorphism (still denoted by $\Phi$) from $\A$ into $\B$. Furthermore, in this case, $\Phi(x)$ is the SOT-limit  of $\{\Phi(pxp)\}_{p\in \paf}$ for any $x \in \A$, and $\norm{\Phi(x)}_\B\leq \norm{x}_\A$ for all $x\in \A^{sa}$.
\end{enumerate}
\end{thm}

Using this result we provide a preliminary structural description of projection disjointness preserving isometries.

\begin{thm} \label{LT9 15/08/16} \label{nL1 15/08/16} 
Suppose $E\subseteq \tma$ and $F\subseteq \nmb$ are symmetric spaces. If $V:E\ra F$ is a projection disjointness preserving isometry, then letting $\Psi(p):=s(V(p))$ for $p\in \paf$, yields a projection mapping that can be extended to a positive linear map (also denoted by $\Psi$) from $\mathcal{F}(\tau)$ into $\B$ such that $\norm{\Psi(x)}_\B=\norm{x}_\A$  and  $\Psi(x^2)=\Psi(x)^2$ for all $x\in \mathcal{F}(\tau)^{sa}$. Furthermore, for any $x\in \mathcal{F}(\tau)$ and $p\in \paf$ with $p \geq s(x)\jn r(x)$, we have 
\begin{enumerate}
\item $V(x)=V(p)\Psi(x)$  
\item $b_p\Psi(x)=\Psi(x)b_p$, where $V(p)=v_{p}b_{p}$ is the polar decomposition of $V(p)$ into a partial isometry $v_p$ and positive operator $b_{p}=|V(p)|$.
\end{enumerate}

\begin{proof}
For $p\in \paf$, let $\Psi(p)=s(V(p))=v_{p}^*v_{p}$. If $p,q\in \paf$ with $pq=0$, then $V(p)^*V(q)=0=V(p)V(q)^*$ and so, as in the proof of Yeadon's Theorem (\cite[Theorem 2]{key-Y81}), we have that that $v_{p}^*v_{q}=0=v_{p}v_{q}^*$. Furthermore, $v_{p}+v_{q}$ is a partial isometry, $|V(p)+V(q)|=b_{p}+b_{q}$ and $V(p)+V(q)=(v_{p}+v_{q})(b_{p}+b_{q})$ is the polar decomposition of $V(p+q)=V(p)+V(q)$. Therefore $v_{p}+v_{q}=v_{p+q}$ and $b_{p}+b_{q}=b_{p+q}$. It follows that 
\[\Psi(p+q)=v_{p+q}^*v_{p+q} =(v_{p}+v_{q})^*(v_{p}+v_{q})=v_{p}^*v_{p}+v_{q}^*v_{q}=\Psi(p)+\Psi(q).\] Using \cite[Exercise 2.3.4]{key-Con1} we have that $\Psi(p)\Psi(q)=0$. Furthermore, if $0\neq p \in \paf$, then $V(p)\neq 0$, since $V$ is injective.  It follows that $\Psi(p)=s(V(p))\neq  0$. Furthermore, by Corollary \ref{CP8 05/08/14}, $V$ has the property that $V(x_n)\rax{\mtop} V(x)$ whenever $\seq{x}$ is a sequence in $\mathcal{F}(\tau)$  such that $x_n \rax{\A} x\in \mathcal{F}(\tau)$ and $s(x_n)\leq s(x)$ for all $n\in \mathbb{N}^+$.  By Theorem \ref{nP extension2}, $\Psi$ can therefore be extended to a positive linear map (also denoted by $\Psi$) from $\mathcal{F}(\tau)$ into $\B$ with the desired properties.

Next we prove (1).  Since $\Psi(p)=s(b_{p})=r(b_{p})=s(v_{p})$, we have that
\begin{eqnarray}
\Psi(p)b_{p}=b_{p}=b_{p}\Psi(p) \qquad \text{and} \qquad v_{p}\Psi(p)=v_{p}. \label{e0.0.1b 15/08/16}
\end{eqnarray}
Suppose $x=q\in \paf$ and $p\in \paf$ with $p\geq q$. Then $p-q\in \paf$ and $q(p-q)=0$. Note that 
$b_{p-q}\Psi(q)=(b_{p-q}\Psi(p-q))\Psi(q)=0$, using (\ref{e0.0.1b 15/08/16}) and the fact that $q(p-q)=0$ implies that $\Psi(q)\Psi(p-q)=0$. Similarly, we have that $v_{p-q}b_{q} =v_{p-q}\Psi(p-q)\Psi(q)b_{q}=0$. Therefore, 
\[V(p)\Psi(q)=v_{q+(p-q)}b_{q+(p-q)}\Psi(q) =(v_{q}+v_{p-q})(b_{q}+b_{p-q})\Psi(q) =v_{q} b_{q}\Psi(q) =V(q).\]
Using the linearity of $V$ and $\Psi$, we therefore have that $V(x)=V(p)\Psi(x)$ for any $x\in \mathcal{G}(\mathcal{A})_f$ and $p\in \paf$ with $p\geq s(x)$. Suppose $x\in \mathcal{F}(\tau)^{sa}$ and $p\in \paf$ with $p\geq s(x)$. As a consequence of the Spectral Theorem, we can find a sequence $\seq{x}$ in $\mathcal{G}(\mathcal{A})_f^{sa}$ such that $x_n \rax{\A} x$ and $s(x_n)\leq s(x)\leq p$ for all $n\ \in \mathbb{N}^+$. By Proposition \ref{P8 05/08/14}, this implies that $x_n \rax{E} x$. Therefore $V(x_n)\rax{F} V(x)$ and $\Psi(x_n)\rax{\B} \Psi(x)$, since $V$ is an isometry and $\Psi$ is linear, and isometric on self-adjoint elements in $\mathcal{F}(\tau)$. Furthermore, since $F$ is a normed $\B$-bimodule, \[\norm{V(p)(\Psi(x_n)-\Psi(x))}_F\leq \norm{V(p)}_F\norm{\Psi(x_n)-\Psi(x)}_\B \ra 0\] and so $V(p)\Psi(x_n) \rax{F} V(p)\Psi(x)$. However, $V(p)\Psi(x_n)=V(x_n)\rax{F} V(x)$. It follows that $V(x)=V(p)\Psi(x)$. Finally, if $x\in \mathcal{F}(\tau)$ and $p\in \paf$ with $p\geq s(x) \jn r(x)$, then $p\geq s(\text{Re}\,x),s(\text{Im}\,x)$ and so $V(x)=V(p)\Psi(x)$ using the linearity of $V$ and $\Psi$.

To prove (2), suppose $x=q$ and $p\in \paf$ with $p\geq q$. Then $b_{p-q}\Psi(q)=b_{p-q}\Psi(p-q)\Psi(q)=0$ and $\Psi(q)b_{p-q}=\Psi(q)\Psi(p-q)b_{p-q}=0$. We therefore have that 
\begin{eqnarray*}
b_{p}\Psi(q)&=&b_{q+(p-q)}\Psi(q) =(b_{q}+b_{p-q})\Psi(q) =b_{q}\Psi(q)\\
& =&\Psi(q)b_{q}=\Psi(q)(b_{q}+b_{p-q}) =\Psi(q)b_{p}.
\end{eqnarray*} 
Noting that for any $p\in \paf$, $b_{p}=v_{p}^*V(p) \in F$ (since $V(p)\in F$, $v_{p}^*\in \B$ and $F$ is a bimodule), we can employ a similar strategy to the one used in (1) to complete the proof. 
\end{proof}
\end{thm}

The previous result allows us to completely describe the structure of projection disjointness preserving isometries in the setting where the initial von Neumann algebra is equipped with a finite trace.

\begin{cor}\label{C2 07/12/18}
Suppose $E\subseteq \tma$ and $F\subseteq \nmb$ are symmetric spaces, and that  $\tau(\id)<\infty$. If $V:E\ra F$ is a projection disjointness preserving isometry, then there exists a Jordan $*$-homomorphism $\Psi$ from $\A$ into $\B$ such that $V(x)=V(\id)\Psi(x)$ for every $x\in \A$.
\end{cor}

For the remainder of this section we will suppose that $(\A,\tau)$ and $(\B,\nu)$ are arbitrary semi-finite von Neumann algebras, $E\subseteq \tma$ and $F\subseteq \nmb$ are symmetric spaces and $U:E\ra F$ is a finiteness preserving positive surjective isometry. It follows from Lemma \ref{L1 C3.1} that $U$ is normal. We will show that there exists a Jordan $*$-isomorphism $\Phi$ from $\A$ onto $\B$ a positive operator $a\in \nmb$ such that
\[U(x)=a\Phi(x)\qquad x\in \A.\]
By Corollary \ref{CnnT3.1 Chilin89}, $U$ is projection disjointness preserving and therefore, by Theorem \ref{nL1 15/08/16}, letting $\Phi(p):=s(U(p))$ for $p\in \paf$ yields a projection mapping which can be extended to a positive linear map (still denoted by $\Phi$) from $\mathcal{F}(\tau)$, which preserves squares of self-adjoint elements. Since $U$ is finiteness preserving and normal, Theorem \ref{T2 04/07/18} can be used to extend $\Phi$ to a normal Jordan $*$-homomorphism (still denoted by $\Phi$) from $\A$ into $\B$. We need to show that $\Phi$ is surjective and define the element $a$ to be used in the representation of $U$. The following lemma will play an important role in both. For $p\in \paf$, we will let $a_p:=U(p)$.

\begin{lem}\label{L1 26/11/18} %new
For any $p\in \paf$, $\Phi(p\A p)=\Phi(p)\B \Phi(p)$.
\begin{proof}
Since $\Phi(pxp)=\Phi(p)\Phi(x)\Phi(p)$ for any $x\in \mathcal{A}$ (see \cite[Exercise 10.5.21]{key-K2}), we have that $\Phi(p\A p)\subseteq\Phi(p)\B \Phi(p)$. Let $y\in \Phi(p)\B \Phi(p)^+$ and define $c=a_p^{1/2}ya_{p}^{1/2}$. Then, since $0\leq y\leq \norm{y}_{\B}\Phi(p)$, repeated application of \cite[Proposition 1(iii)]{key-Dodds14} yields
\begin{eqnarray}
0 \leq c \leq a_p^{1/2}\norm{y}_{\B}\Phi(p)a_p^{1/2} = \norm{y}_{\B}a_p,  \label{e4.0 C3.1b}
\end{eqnarray}
using the fact that $\Phi(p)=s(a_p)= s(a_p^{1/2})$. Since $F$ is symmetric (and hence absolutely solid) and $\norm{y} a_p=\norm{y}U(p)\in F$, it follows that $c\in F$. By Lemma \ref{L1 C3.1}, $U^{-1}$ is positive and therefore $0\leq U^{-1}(c)\leq \norm{y}p$. It follows that $U^{-1}(c)\in p\A p$. By Theorem \ref{T2 04/07/18}, there exists a $w_p\in \nmb$ such that $U(p)^{1/2}w_p=\Phi(p)=w_pU(p)^{1/2}$ and $\Phi(x)=w_pU(x)w_p$. Since $a_p=U(p)$, it follows that 
\[\Phi(U^{-1}(c))=w_pU(U^{-1}(c))w_p=w_p(a_p^{1/2}ya_p^{1/2})w_p=\Phi(p)y\Phi(p)=y.\]
Since elements in $\Phi(p)\B \Phi(p)$ can be written as finite linear combinations of elements $\Phi(p)\B \Phi(p)^+$, we have that $\Phi(p)\B \Phi(p)\subseteq \Phi(p\A p)$.
\end{proof}
\end{lem}

Next we define $a$. Let $a_p=\int_0^{\infty}\lambda de^{a_p}_\lambda$ denote the spectral representation of $a_p$. We start by showing that for a fixed $\lambda \geq 0$, $\{e^{a_p}(\lambda,\infty)\}_{p\in \paf}$ is an increasing net, where $e^{a_p}(\lambda,\infty)=\id-e_\lambda^{a_p}$. Suppose $q\in \paf$ with $q\geq p$. Note that $e^{a_p}(\lambda,\infty)\leq s(a_p)=\Phi(p)\leq \Phi(q)$ and so, by Lemma \ref{L1 26/11/18}, there exists an $x\in q\A q$ such that $e^{a_p}(\lambda,\infty)=\Phi(x)$. It follows by Theorem \ref{nL1 15/08/16}(2) that $a_q e^{a_p}(\lambda,\infty)= e^{a_p}(\lambda,\infty)a_q$ and therefore $e^{a_q}(\lambda,\infty)e^{a_p}(\lambda,\infty)=e^{a_p}(\lambda,\infty)e^{a_q}(\lambda,\infty)$. Since $U$ is positive, we also have that $a_p=U(p)\leq U(q)=a_{q}$. Therefore $e^{a_p}(\lambda, \infty)\leq  e^{a_{q}}(\lambda, \infty)$ for all $\lambda\geq 0$. By \cite[Proposition 2.5.6]{key-K1}, $\{e^{a_p}(\lambda,\infty)\}_{p\in \paf}$ converges in the strong operator topology. Define $e^a(\lambda,\infty):=\sotlimx{p\in \paf}\,\,e^{a_p}(\lambda,\infty)$ and $e_\lambda^a=\id-e^a(\lambda,\infty)$. One can show that $\{e^a_\lambda\}_{\lambda \geq 0}$ is a resolution of the identity and, by \cite[Lemma 5.6.9]{key-K1},  letting 
\[a=\int_0^{\infty}\lambda de^a_\lambda\] yields a closed and densely defined positive operator. Furthermore $a_p=U(p)\in F\subseteq \nmb$ and so $e^{a_p}_\lambda\in \B$ for each $\lambda \geq 0$. Since $\B$ is closed in the strong operator topology, it follows that $e^a_\lambda \in \B$ for each $\lambda \geq 0$ and therefore $a$ is affiliated with $\B$. Before discussing the relationship between $a$ and $\Phi$, which will enable us to show that $a\in \nmb$, we include a result that we will need. It is likely that this is a known result, but since the authors were unable to find an appropriate reference we also include a short proof.

\begin{prop} \label{P2 19/03/14}
Let $x$ be a closed, densely defined self-adjoint operator on $H$ with spectral representation  $x=\int_{-\infty}^{\infty}\lambda de^x_\lambda$. If $p$ is a projection such that $px= xp$, then   $px=\int_{-\infty}^{\infty} \lambda d(p e^x_\lambda)$ (i.e. $\{p e^x_\lambda \}_\lambda$ is the resolution of the identity for $px$).
\begin{proof}
Let $\{e_\lambda^x\}_{\lambda \in \mathbb{R}}$ denote the resolution of the identity for $x$. For each $n\in \n$, put $f^x_n=e^x_n - e^x_{-n}$. Then for each $n$ and each $\xi\in f^x_n(H)$, $x\xi=\int_{-n}^n\!\lambda  d e^x_\lambda \xi$ (\cite[Lemma 5.6.7]{key-K1}). Since $px= xp$, $p$ commutes with $e_\lambda^x$ for each $\lambda\in \br$, by \cite[Theorem 1.5.12]{key-dP}, and so $e_\lambda:=e_\lambda^x p$ is a projection for each $\lambda$. It is easily checked that $\{e_\lambda\restriction_{p(H)}\}_{\lambda \in \mathbb{R}}$ is a resolution of the identity on the Hilbert space $pHp$. It follows, using the fact that the integral is a limit of linear combinations of disjoint spectral projections commuting with $p$, that $px\xi =xp\xi= (\int_{-n}^n \lambda de^x_\lambda)p\xi=\int_{-n}^{n} \lambda d(e^x_\lambda p)\xi$ for $n\in\n$ and  $\xi\in f_n(H)$, where $f_n= e_n-e_{-n}$. Since $\unionx{n=1}{\infty} f_n(H)$ is a core for $xp$, the result follows by \cite[Theorem 5.6.12]{key-K1}.
\end{proof}
\end{prop}

We return now to discussing the relationship between $a$ and $\Phi$.

\begin{lem} \label{L7.2 C3.1}\label{LeL7.3 C3.1}
If $p\in \paf$, then $a\Phi(p)=a_p=\Phi(p)a$.
\begin{proof}
We start by showing that $e^a_\lambda\Phi(p)=e^{a_p}_\lambda=\Phi(p)e^a_\lambda$ for $\lambda>0$ and $p\in \paf$. Let $q\in \paf$ with $q\geq p$. Then, using the definition of $a_p$ and applying Theorem \ref{nL1 15/08/16}, we obtain
\begin{eqnarray}
a_p=U(p) =U(q)\Phi(p)= a_{q}\Phi(p) \label{e4.11 C3.1}
\end{eqnarray}
Furthermore, $\Phi(p)$ is a projection and $a_{q}\Phi(p)=\Phi(p)a_{q}$, by Theorem \ref{nL1 15/08/16}(2). Using Proposition \ref{P2 19/03/14} and (\ref{e4.11 C3.1}) it follows that $\{e^{a_{q}}_\lambda \Phi(p)\}_\lambda$ is the resolution of the identity for $a_q\Phi(p)=a_p$, i.e. $e^{a_p}_\lambda=e^{a_{q}}_\lambda \Phi(p)$ for every $\lambda \geq 0$. Furthermore,  $e^{a_{q}}_\lambda \Phi(p) \sotc e^a_\lambda \Phi(p)$ as $q \uparrow \id$. Therefore, $e^{a_p}_\lambda= e^a_\lambda \Phi(p)$. Since $a_q\Phi(p)=\Phi(p)a_q$, we have that $e^{a_q}_\lambda\Phi(p)=\Phi(p)e^{a_q}_\lambda$ and therefore appropriate adjustments to the last few lines yields $e^{a_p}_\lambda=\Phi(p)e^a_\lambda$. Combining this with what was shown earlier we obtain $e^a_\lambda \Phi(p)=e^{a_p}_\lambda=\Phi(p)e^a_\lambda$. Therefore, using a similar approximation argument to the one employed at the end of Proposition \ref{P2 19/03/14}, we obtain 
\[ a_p=\int_0^{\infty}\lambda de^{a_p}_\lambda=\int_0^{\infty}\lambda d(e^a_\lambda\Phi(p)) = (\int_0^{\infty}\lambda de^a(\lambda))\Phi(p)=a\Phi(p). \]
Similarly, $a_p=\Phi(p)a$.
\end{proof}
\end{lem}

Since $\Phi(p)$ is defined everywhere, $\D(\Phi(p)a)=\{\eta\in \D(a):a\eta \in \D(\Phi(p))\}=\D(a)$. It follows that $\D(a)=\D(\Phi(p)a)=\D(a_p)$ and therefore $\D(a)$ is $\nu$-dense, since $a_p=U(p)\in F \subseteq \nmb$. Thus $a\in \nmb$, since we have already shown that $a$ is a closed densely defined operator affiliated with $\B$.

\begin{lem} \label{LP8.5 C3.1} 
 If $x\in \A \cap E$, then $U(x)=a\Phi(x)$.
\begin{proof}
Suppose $x\in \mathcal{F}(\tau)^{sa}$ and let $p=s(x)$. Then $r(\Phi(x))\leq \Phi(p)$, by \cite[Lemma 3.5]{key-dJ18b}. Using  Theorem \ref{nL1 15/08/16} and Lemma \ref{LeL7.3 C3.1}, we therefore have $U(x)=a_p\Phi(x)=a\Phi(p)\Phi(x)=a\Phi(x)$. Next, suppose that $x\in \A^+\cap E$. By \cite[Proposition 1(vii)]{key-Dodds14} there exists an increasing net $\net{x}$ in $\mathcal{F}(\tau)^+$ such that $x_\lambda \uparrow x$. Then using the normality of $U$ and $\Phi$ we have that $U(x_\lambda)\uparrow U(x)$ and $\Phi(x_\lambda)\uparrow \Phi(x)$. Therefore $U(x_\lambda) \lmtopc U(x)$ and $\Phi(x_\lambda) \lmtopc \Phi(x)$, by \cite[Proposition 2(v)]{key-Dodds14}. It follows that $a\Phi(x_\lambda) \lmtopc a\Phi(x)$ (see \cite[p.211]{key-Dodds14}). Since $a\Phi(x_\lambda)=U(x_\lambda)$ for each $\lambda$ and the local measure topology is Hausdorff (\cite[Proposition 2.7.4]{key-dP}), we have that $U(x)=a\Phi(x)$.
\end{proof}
\end{lem}

\begin{lem}\label{L2 26/11/18} 
$\Phi$ is a Jordan $*$-isomorphism from $\A$ onto $\B$
\begin{proof}
Assume that $\id -\Phi(\id)\neq 0$. Since $(\B,\nu)$ is semi-finite, there exists a $q\in \px{B}$ such that $0<q\leq \id - \Phi(\id)$ and $\nu(q)<\infty$. This implies that $q\in F$ and hence there exists an $x\in E$ such that $U(x)=q$, since $U$ is surjective. By \cite[Proposition 1(vii)]{key-Dodds14}, there exists $\net{x}$ in $\mathcal{F}(\tau)$ such that $x_\lambda \uparrow  x$. Then, using Lemma \ref{LP8.5 C3.1} and the normality of $U$, we obtain $a\Phi(x_\lambda)=U(x_\lambda)\uparrow U(x)=q$. Therefore $a\Phi(x_\lambda) \lmtopc q$. However, we also have that $a\Phi(x_\lambda)=a\Phi(x_\lambda)\Phi(\id)\lmtopc q\Phi(\id)$, by \cite[Exercise 10.5.22]{key-K2}, Lemma \ref{LP8.5 C3.1} and \cite[p. 211]{key-Dodds14}. It follows that $q=q\Phi(\id)$. However, since $q\leq \id - \Phi(\id)$, we have that $q(\id-\Phi(\id))=q$, and so $q=q\Phi(\id)=\Bigl(q(\id-\Phi(\id))\Bigr)\Phi(\id)=0$. This is a contradiction and so $\Phi$ is unital.

Noting that \cite[Theorem 4.5]{key-dJ18b} is employed in the proof of \cite[Theorem 5.1]{key-dJ18b} and considering \cite[Remark 4.6]{key-dJ18b}, it follows that $\Phi$ is isometric on $\A^{sa}$, since $\Phi(p)=s(U(p))=0$ if and only if $p=0$. By Lemma \ref{L1 26/11/18}, $\Phi(p)\B \Phi(p)= \Phi(p\A p)\subseteq \Phi(\A)$ for every $p\in \paf$ and therefore $\Phi$ is a Jordan $*$-isomorphism from $\A$ onto $\B$, by \cite[Proposition 6.2]{key-dJ18b}.  
\end{proof}
\end{lem}

We have therefore obtained the following result.

\begin{thm} \label{T1 23/11/18} 
Suppose $(\A,\tau)$ and $(\B,\nu)$ are semi-finite von Neumann algebras, $E\subseteq \tma$ and $F\subseteq \nmb$ are symmetric spaces and $U:E\ra F$ is a positive surjective isometry. If $U$ is finiteness preserving (in particular if $\nu(\id)<\infty$), then there exists a positive operator $a\in \nmb$ and a Jordan $*$-isomorphism $\Phi$ of $\A$ onto $\B$ such that $U(x)=a\Phi(x)$ for all $x\in \A\cap E$.
\end{thm}

\section{The structure of projection disjointness preserving isometries} \label{S5}

In the previous section we showed that under certain conditions the structure of a positive surjective isometry can be described in terms of a positive operator and Jordan $*$-isomorphism. We will use this result to show that we can obtain a similar representation for a surjective isometry, which is not necessarily positive, if it is projection disjointness preserving. Throughout this section we will suppose that $E\subseteq \tma$ is a strongly symmetric space with absolutely continuous norm,  $F\subseteq \nmb$ is a symmetric space and $V:E\ra F$ is a projection disjointness and finiteness preserving surjective isometry. The idea of the proof, inspired by \cite[$\S 5$]{key-Chilin89}, is to use the isometry $V$ to construct a unitary operator $v$ such that $v^*V(\cdot)$ yields a positive surjective isometry and whose structure can therefore be described by the results of the previous section.

By Theorem \ref{nL1 15/08/16},  letting $\Psi(p):=s(V(p))$ for $p\in \paf$, yields a projection mapping that can be extended to a positive linear map (also denoted by $\Psi$) from $\mathcal{F}(\tau)$ into $\B$ with Jordan $*$-homomorphism-like properties (i.e. $\Psi$ is positive, $\norm{\Psi(x)}_\B=\norm{x}_\A$  and  $\Psi(x^2)=\Psi(x)^2$ for all $x\in \mathcal{F}(\tau)^{sa}$). As in Theorem \ref{nL1 15/08/16}, we will, for each $p\in \paf$, write  $V(p)=v_pb_p$ for the the polar decomposition of $V(p)$.

\begin{lem} \label{L3 15/08/16} $\{v_{p}\}_{p\in \paf}$ converges in the strong operator topology to a unitary operator $v\in \B$ and $v\Psi(p) =v_{p}$ for all $p\in \paf$.
\begin{proof}
We start by noting that if $p,q\in \paf$ are such that $0<q\leq p$, then \begin{eqnarray}
v_{q}=v_{p}\Psi(q).  \label{eL2 15/08/16}
\end{eqnarray} To show this, note that if $p=q$, then (\ref{eL2 15/08/16}) holds using (\ref{e0.0.1b 15/08/16}). If $p>q$, then $0\neq p-q \in \paf$ and $q(p-q)=0$. Therefore, $v_{q}+v_{p-q}=v_{p}$ (see \cite[Proposition B.1.32(5)]{key-dJ17}). It follows that $ v_{p}\Psi(q)=(v_{q}+v_{p-q})\Psi(q) =v_{q} +v_{p-q}\Psi(p-q)\Psi(q) =v_{q}$, using (\ref{e0.0.1b 15/08/16}) and the fact that $(p-q)q=0$ implies that $\Psi(p-q)\Psi(q)=0$.

Next, we show that $\{v_{(p)}\}_{p\in \paf}$ is SOT-convergent to a partial isometry. Let $\eta\in K$ (where $\B\subseteq B(K)$) and suppose $\epsilon>0$. Since $\{\Psi(p)\}_{p\in \paf}$ is an increasing net of projections, it converges in the strong operator topology to a projection $y\in \pb$. It follows that there exists a $p_\epsilon \in \paf$ such that $p,q\in \paf$ with $p,q\geq p_\epsilon$ implies that $\norm{(\Psi(p)-\Psi(q))\eta}<\epsilon$. Let $p,q\in \paf$ with $p,q\geq p_\epsilon$. Since $\paf$ is a directed set, there exists an $r\in \paf$ with $r\geq p,q$. Using (\ref{eL2 15/08/16}), we then have 
\[\norm{(v_{p}-v_{q})\eta}=\norm{v_{r}(\Psi(p)-\Psi(q))\eta} \leq\norm{v_{r}}_\B\norm{(\Psi(p)-\Psi(q))\eta}<\epsilon. \]
Therefore $\{v_{p}(\eta)\}_{p\in \paf}$ is Cauchy in $K$. Since this holds for every $\eta\in K$, we have that $\{v_{p}\}_{p\in \paf}$ is SOT-Cauchy. Furthermore, $\{v_{p}\}_{p\in \paf}$ is contained in the unit ball of $\mathcal{B}(K)$ and so $v_{p}\rax{SOT}v$ for some $v\in \bhx{K}$, since norm-closed balls in $\bhx{K}$ are SOT-complete by \cite[Proposition 2.5.11]{key-K1}. Since $\B$ is SOT-closed, $v\in \B$. Furthermore, for any $q\in \paf$ with $q\geq p$, we have $v_{p}=v_{q}\Psi(p) \rax{SOT}v\Psi(p)$ as $q\uparrowx{q\in \paf} \id$ using (\ref{eL2 15/08/16}) and the fact that multiplication is separately continuous in the strong operator topology. It follows that $v_{p}=v\Psi(p)$. We show that $v$ is a partial isometry and $s(v)=y$. Note that $v_{p} \rax{SOT} v$ implies that $v_{p} \rax{WOT} v$ since the WOT is coarser than the SOT. Therefore $v_{p}^* \rax{WOT} v^*$ (see \cite[Exercise 5.7.1]{key-K1}) and so $v_{p}^*v_{p} \rax{WOT} v^*v$. Furthermore, $v_{p}^*v_{p}=|v_{p}|=\Psi(p)\rax{SOT} y$ and so $v_{p}^*v_{p}\rax{WOT} y$. It follows from the uniqueness of weak operator topology limits, this implies that $ y=v^*v$. Therefore $v$ is a partial isometry (see \cite[Proposition 6.1.1]{key-K2}) and $s(v)=y$.

We show that $y=\id$ and hence that $v$ is unitary. Suppose $x\in \mathcal{F}(\tau)$. For $p\in \paf$ with $p\geq s(x) \jn r(x)$ we have that $px=x=xp$ and hence $\Psi(x)=\Psi(xp)=\Psi(x)\Psi(p)$ (see \cite[Proposition 2.3]{key-dJ18b}). Therefore,  \[\Psi(x)y=\Psi(x) \sotlimx{p\in \paf} \Psi(p) =\sotlimx{p\in \paf:p\geq s(x)\jn r(x)}[\Psi(x)\Psi(p)] =\Psi(x).\] It follows that if $p\geq s(x) \jn r(x)$, then $V(x)=b_{p} v_{p} \Psi(x) =b_{p} v_{p} \Psi(x)y =V(x)y$, using Theorem \ref{nL1 15/08/16}. Assume that $\id -y \neq 0$. Since, $(\B,\nu)$ is semi-finite, there exists a $q\in \px{B}$ such that $0<q\leq \id - y$ and $\nu(q)<\infty$. This implies that $q\in F$ and hence there exists an $x\in E$ such that $V(x)=q$, since $V$ is surjective. $E$ has absolutely continuous norm and therefore $\mathcal{F}(\tau)$ is dense in $E$ (see \cite[p.241]{key-Dodds14}). Let $\seq{x}$ be a sequence in $\mathcal{F}(\tau)$  such that $x_n \rax{E} x$. Then $V(x_n) \rax{F} V(x)=q$. However $V(x_n)=V(x_n)y \rax{F} V(x)y=qy$ and so $q=qy =0$, since $q\leq \id - y$. This is a contradiction and so $y=\id$. 
\end{proof}
\end{lem}

\begin{lem} \label{L8 15/08/16}
The map $U:E \ra F$ defined by $U(x)=v^*V(x)$ is a positive surjective isometry. 
\begin{proof}
Since $v^*$ is a unitary operator, it is easily checked that $U$ is a surjective isometry. To see that $U$ is positive note that if $x\in \mathcal{F}(\tau)^+$ and $p=s(x)$, then $p\in \paf$ and $V(x)=v_{p}b_{p}\Psi(x)=v\Psi(p)b_{p}\Psi(x) =vb_{p}\Psi(x)$, by Theorem \ref{nL1 15/08/16}, Lemma \ref{L3 15/08/16} and (\ref{e0.0.1b 15/08/16}). It follows that $v^*V(x)=b_{p}\Psi(x) =b_{p}^{1/2}\Psi(x)b_{p}^{1/2} \geq 0$ using Theorem \ref{nL1 15/08/16}, \cite[Proposition 2.2.22]{key-dP} (with $f(t):=t^{1/2}$),  \cite[Proposition 1(iii)]{key-Dodds14} and the fact that $\Psi$ is positive.
Suppose $x\in E^+$. Since $E$ has absolutely continuous norm, there exists a sequence $\seq{x}$ in $\mathcal{F}(\tau)^+$ such that $x_n \rax{E} x$. As $U$ is an isometry,  $v^*V(x_n)=U(x_n)\rax{F} U(x)$. We have that $v^*V(x_n)\geq 0$ for all $n\in \mathbb{N}^+$ and therefore $U(x)\geq 0$, since $F^+$ is closed by \cite[Corollary 12(i)]{key-Dodds14}.  
\end{proof}
\end{lem}

\begin{thm} \label{T9 15/08/16}
Suppose $E$ is a strongly symmetric space with absolutely continuous norm and $F$ is a symmetric space. If $V:E\ra F$ is a projection disjointness and finiteness preserving surjective isometry, then there exists a unitary operator $v$, a positive operator $a$ affiliated with the centre of $\B$ and a Jordan $*$-isomorphism $\Phi$ from $\A$ onto $\B$ such that $V(x)=va\Phi(x)$ for all $x\in \A \cap E$. 
\begin{proof}
In order to apply Theorem \ref{T1 23/11/18} to describe the structure of $U$ as defined by the previous lemma we need to show that $U$ is finiteness preserving. To this end, suppose that $p\in \paf$. Then \[U(p)=v^*V(p)=v^*v_{p}b_{p}=v^*v\Psi(p)b_{p} =b_{p},\]
by Theorem \ref{nL1 15/08/16}, Lemma \ref{L3 15/08/16} and (\ref{e0.0.1b 15/08/16}). It follows from the above and the finiteness preserving assumption on $V$ that $\nu(s(U(p)))=\nu(s(b_{p}))=\nu(s(V(p)))<\infty$.  By Theorem \ref{T1 23/11/18},  there exists a positive operator $a\in \nmb$ and a Jordan $*$-isomorphism $\Phi$ from $\A$ onto $\B$ such that $U(x)=a\Phi(x)$ for all $x\in \A \cap E$ and so $V(x)=va\Phi(x)$ for all $x\in \A \cap E$.
\end{proof}
\end{thm}

\begin{rem}
We demonstrate briefly that $\Phi$ (obtained in the theorem above) is the unique normal extension of $\Psi:\mathcal{F}(\tau)\ra \B$ (as obtained earlier in this section by extending the map $\Psi(p):=s(V(p))$ for $p\in \paf$) and that $b_{p}=a_{p}$ for every $p\in \paf$, where the $a_{p}$'s are the positive operators used to construct $a$ as in $\S \ref{S4}$. Recall that $a_p=U(p)=\int_0^\infty \lambda de^{a_p}_\lambda$, $e^a_\lambda=\sotlimx{p\in \paf}\,\,e^{a_p}_\lambda$ and $a=\int_0^\infty \lambda de^a(\lambda)$. However, $b_{p}=v^*V(p)=U(p)$ and so $b_{p}=a_p$ for every $p\in \paf$. To demonstrate the relationship between $\Phi$ and $\Psi$, recall that $\Phi$ is obtained using Theorem \ref{nP extension2} and as such $\Phi(p)=s(U(p))=s(v^*V(p))=s(V(p))=\Psi(p)$ for every $p\in \paf$, since $v^*$ is unitary. 
\end{rem}

\section*{Acknowledgments}

The greater part of this research was conducted during the first author's doctoral studies at the University of Cape Town. The first author would like to thank his Ph.D. supervisor, Dr Robert Martin, for his input and guidance. Furthermore, the support of the DST-NRF Centre of Excellence in Mathematical and Statistical Sciences (CoE-MaSS) towards this research is hereby acknowledged. Opinions expressed and conclusions arrived at, are those of the authors and are not necessarily attributed to the CoE.

\bibliographystyle{amsplain}

\end{document}